\documentclass[graybox]{svmult}

\usepackage{mathptmx}       %
\usepackage{helvet}         %
\usepackage{courier}        %
\usepackage{type1cm}        %
\usepackage{makeidx}         %
\usepackage{graphicx}        %
\usepackage{multicol}        %
\usepackage[bottom]{footmisc}%

\usepackage{epsfig}
\usepackage{graphicx}
\usepackage[latin1]{inputenc}

\usepackage{amsfonts,bbm,amssymb,epsf,amsmath,%
latexsym,amscd,enumerate,supertabular,color,url}

\usepackage{subfigure,psfrag}

\newcommand{\p}{\partial}

\usepackage{algorithm}
\usepackage{algorithmic}

\newcommand{\mpr}[1]{\relax}
\newcommand{\mprup}[1]{\relax}
\newcommand{\mprdown}[1]{\relax}

\newcommand{\pf}{{\noindent \it Proof.} }
\newcommand{\eop}{{\qed}\par\bigskip}
\newcommand{\meop}{\quad{\qed}}

\def\sp{\,:\;}
\def\ie{{\sl \thinspace i.e.},\ }

\def\tri{\triangle}

\newcommand{\RR}{\mathbb{R}}

\newcommand{\CN}{\mathcal{N}}

\newcommand{\sps}{\mathbb{S}^{1,2}_{5,0}(\tri)}

\makeindex             %

\begin{document}

\title*{Approximation by $C^1$ Splines
on Piecewise Conic Domains}
\author{Oleg Davydov and Wee Ping Yeo}
\institute{Oleg Davydov \at Department of
Mathematics, University of Giessen,
    Department of Mathematics,
    Arndtstrasse 2,
    35392 Giessen,
    Germany, \email{oleg.davydov@math.uni-giessen.de}
\and Wee Ping Yeo \at Faculty of Science,
Universiti Brunei Darussalam,
BE1410, Brunei Darussalam, \email{weeping.yeo@ubd.edu.bn}}
\maketitle

\newcommand{\abstrtext}{We develop a Hermite interpolation scheme and prove error bounds for $C^1$
bivariate piecewise polynomial spaces of Argyris type vanishing on the boundary of
curved domains enclosed by piecewise conics. }

\abstract*{\abstrtext}

\abstract{\abstrtext}

\section{Introduction}
\label{intro}
Spaces of piecewise polynomials defined on domains bounded by piecewise algebraic curves and vanishing on
parts of the boundary can be used in the Finite Element Method as an alternative to the classical mapped
curved elements \cite{DKS,DSaC1}. Since implicit algebraic curves and surfaces provide a well-known
modeling tool in CAGD \cite{ImplSurf}, these methods are inherently isogeometric in the sense of \cite{HCB05}.
Moreover, this approach does not suffer from the usual difficulties
 of building a globally $C^1$ or smoother space of functions on curved domains
(see \cite[Section 4.7]{BrennerScott})
shared by  the classical curved finite elements and the
B-spline-based isogeometric analysis.

In particular, a space of $C^1$ piecewise polynomials on domains enclosed by piecewise conic sections has been studied in
\cite{DSaC1} and applied to the numerical solution of fully nonlinear elliptic equations. These piecewise polynomials are
quintic on the interior triangles of a triangulation of the domain, and sextics on the boundary triangles
(pie-shaped triangles with one side represented by a conic section as well as
those triangles that share with them an interior edge with one endpoint on the boundary) and
generalize the well-know Argyris finite element. Although local bases for these spaces have been constructed in
\cite{DSaC1} and numerical examples demonstrated the convergence orders expected from a piecewise
quintic finite element, no error bounds have been provided.

In this paper we study the approximation properties of the spaces introduced in
\cite{DSaC1}. We define a Hermite-type interpolation operator and prove an error bound that shows the
convergence order $\mathcal{O}(h^6)$ of the residual in $L_2$-norm, and order $\mathcal{O}(h^{6-k})$ in Sobolev
spaces $H^k(\Omega)$. This extends the techniques used in \cite{DKS} for $C^0$ splines to
Hermite interpolation.

The paper is organized as follows. We introduce in Section 2 the spaces $\mathbb{S}^{1,2}_{d,0}(\tri)$ of $C^1$ piecewise
polynomials on
domains bounded by a number of conic sections, with homogeneous boundary conditions, define  in Section 3 our interpolation
operator in the case $d=5$,  and investigate in Section 4 its
approximation error  for functions in Sobolev spaces $H^m(\Omega)$, $m=5,6$, vanishing on the boundary.

\section{$C^1$ piecewise polynomials on  piecewise conic domains\label{spacesC1}}\mprup{spaces}

We make the same assumptions on the domain and its triangulation as in \cite{DKS,DSaC1}, as outlined below.

Let $\Omega\subset\RR^2$ be a bounded curvilinear polygonal domain with
$\Gamma=\partial \Omega=\bigcup_{j=1}^n\overline{\Gamma}_j$,
where each $\Gamma_j$ is an open arc of an algebraic curve of
at most second order
(\ie either a straight line or a conic). For simplicity we assume that
$\Omega$ is simply connected, so that its boundary $\Gamma$ is a closed curve without self-intersections.
Let $Z=\{z_1,\ldots,z_n\}$ be the set of the endpoints of all
arcs numbered counter-clockwise such that $z_j,z_{j+1}$ are the
endpoints of $\Gamma_j$, $j=1,\ldots,n$, with $z_{j+n}=z_j$.
Furthermore, for each $j$ we denote by $\omega_j$ the internal
angle between the tangents $\tau^+_j$ and $\tau^-_j$ to $\Gamma_j$
and $\Gamma_{j-1}$, respectively, at $z_j$. We assume that $\omega_j\in(0,2\pi)$ for all $j$. Hence
$\Omega$ is a Lipschitz domain.

\begin{figure}[htbp!]
\centering
\begin{tabular}{c}
\includegraphics[height=0.3\textwidth]{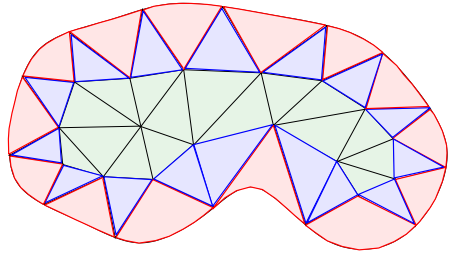}
\end{tabular}
\caption{
A triangulation of a curved domain with ordinary triangles (green),
pie-shaped triangles (pink) and buffer triangles (blue).
}\label{curved_tri}
\end{figure}

Let $\tri$ be a \emph{triangulation} of $\Omega$, \ie a
subdivision of $\Omega$ into triangles, where
each triangle $T\in\tri$ has at most one edge replaced with a curved segment of
the boundary $\partial\Omega$, and the intersection of any pair of the triangles
is either a common vertex or a common (straight) edge if it is non-empty.
The triangles with a curved edge are said to be \emph{pie-shaped}.
Any triangle $T\in\tri$ that shares at least one edge with a pie-shaped triangle
is called a \emph{buffer} triangle, and the remaining triangles are
\emph{ordinary}. We denote by $\tri_0$, $\tri_B$ and $\tri_P$ the sets of all
ordinary, buffer and pie-shaped triangles of $\tri$, respectively, such that
$\tri=\tri_0\cup\tri_{B} \cup \tri_P$
is a disjoint union, see Figure~\ref{curved_tri}.
Let $V,E,V_I,E_I,V_\partial,E_\partial$ denote the set of all vertices,
all edges, interior vertices, interior edges, boundary vertices and boundary  edges,
respectively.

For each  $j=1,\ldots,n$, let $q_j\in\mathbb{P}_2$ be a polynomial such that
$\Gamma_j\subset\{x\in\RR^2\sp q_j(x)=0\}$, where $\mathbb{P}_d$ denotes the space of all bivariate polynomials of total degree
at most $d$. By changing the sign of $q_j$ if needed, we ensure that
$q_j(x)$ is positive for points in $\Omega$ near the boundary segment $\Gamma_j$.
For simplicity we assume in this paper that all boundary segments $\Gamma_j$
are curved. Hence each
$q_j$ is an irreducible quadratic polynomial and
\begin{equation}\label{qgrad}
\nabla q_j(x)\ne 0\quad\text{if}\quad x\in \Gamma_j.
\end{equation}

We assume that $\tri$ satisfies the following conditions:
\begin{itemize}
\item[(A)]\quad $Z=\{z_1,\ldots,z_n\}\subset V_\partial$.
\item[(B)]\quad No interior edge has both endpoints on the boundary.
\item[(C)]\quad No pair of pie-shaped triangles shares an edge.
\item[(D)]\quad Every $T\in\tri_P$ is star-shaped with respect to its interior vertex $v$.
\item[(E)]\quad  For any  $T\in\tri_P$ with its curved side on $\Gamma_j$,
$q_j(z)>0$ for all $z\in T\setminus \Gamma_j$.
\item[(F)]\quad No pair of buffer triangles shares an edge.
\end{itemize}
It can be easily seen that (B) and (C) are achievable by a slight modification of
a given triangulation, while (D) and (E) hold for sufficiently fine triangulations.
The assumption (F) is made
for the sake of simplicity of the analysis. Note that the triangulation shown in Figure~\ref{curved_tri}
does not satisfy (F).

\begin{figure}
\centering
\begin{tabular}{c}
\psfrag{v1}{$v_1$}
\psfrag{v2}{$v_2$}
\psfrag{v3}{$v_3$}
\psfrag{B}{$B_T$}
\psfrag{T}{$T^*$}
\includegraphics[height=0.35\textwidth]{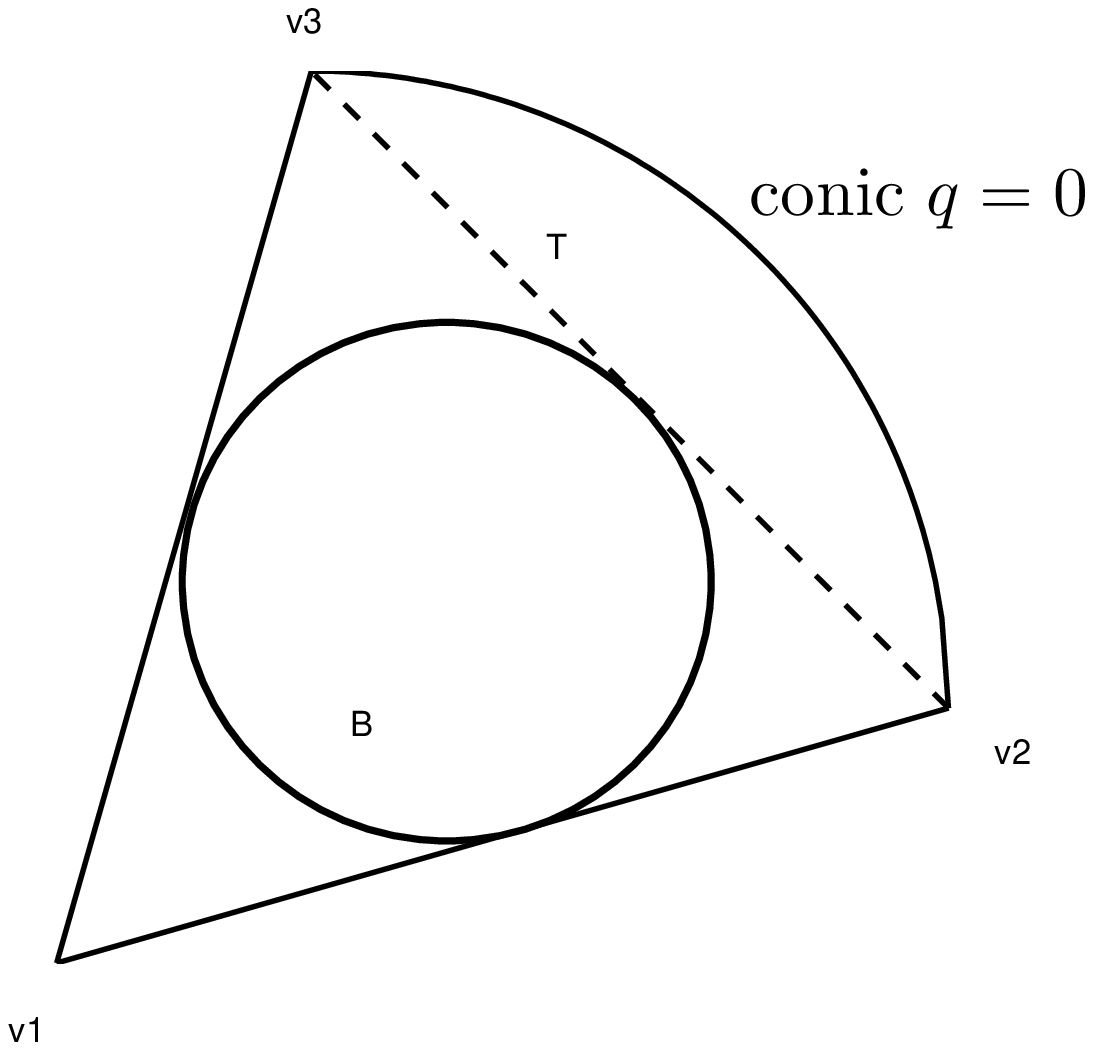}
\qquad
\psfrag{v1}{$v_1$}
\psfrag{v2}{$v_2$}
\psfrag{v3}{$v_3$}
\psfrag{B}{$B_T$}
\psfrag{T2}{$T^*$}
\includegraphics[height=0.35\textwidth]{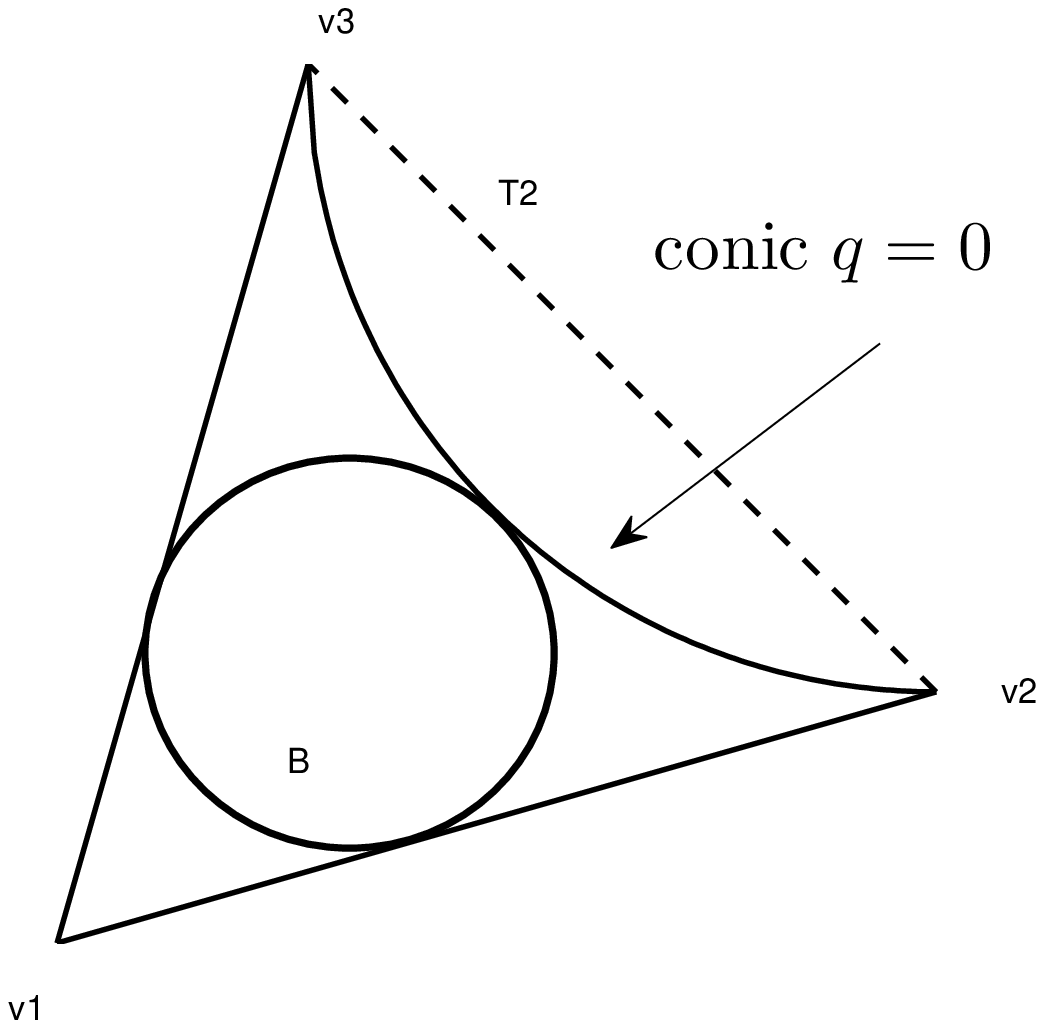}
\end{tabular}
\caption{
A pie-shaped triangle with a curved edge and the
associated triangle $T^\ast$ with straight sides and vertices $v_1,v_2,v_3$.
The curved edge can be either outside (left) or inside $T^\ast$ (right).
}\label{pieshapedT1}
\end{figure}

For any $T\in\tri$, let $h_T$ denote the diameter of $T$, and let
$\rho_T$ be the radius of the disk $B_T$ inscribed in $T$ if $T\in\tri_0\cup\tri_B$ or in $T\cap T^*$ if
$T\in\tri_P$, where $T^\ast$ denotes the triangle obtained by joining the boundary vertices of $T$
by a straight line, see Figure~\ref{pieshapedT1}. Note that every triangle $T\in\tri$ is star-shaped
with respect to $B_T$. %
In particular, for
$T\in\tri_P$ this follows from Condition (D) and the fact that the conics do not possess inflection points.

We define the \emph{shape regularity constant} of $\tri$ by
\begin{equation}\label{shape}
R=\max_{T\in\tri}\frac{h_T}{\rho_T}.
\end{equation}

For any $d\ge1$ we set
\begin{align*}
\mathbb{S}^1_d(\tri)&:=\{s\in C^1(\Omega)\sp s|_T\in\mathbb{P}_{d},\;T\in\tri_0,\;\text{and}\;
s|_T\in\mathbb{P}_{d+1},\;T\in\tri_P\cup\tri_B\},
				 \\
\mathbb{S}^{1,2}_{d,I}(\tri)&:=\{s\in  \mathbb{S}^1_d(\tri)\sp s \text{ is twice differentiable at any }v\in V_I\}
				 ,\\
\mathbb{S}^{1,2}_{d,0}(\tri)&:=\{s\in \mathbb{S}^{1,2}_{d,I}(\tri)\sp s|_\Gamma = 0\}. %
\end{align*}

We refer to \cite{DSaC1} for the construction of a local basis for the space $\mathbb{S}^{1,2}_{5,0}(\tri)$ and its applications in
the Finite Element Method.

Our goal is to obtain an error bound for the approximation of functions vanishing on the
boundary by splines in $\mathbb{S}^{1,2}_{5,0}(\tri)$. This is done through the construction of 
an interpolation operator of Hermite type.
Note that a method of \emph{stable splitting} was employed in \cite{D07,DSa12,DSa13} to estimate the
approximation power of $C^1$ splines vanishing on the boundary of a polygonal domain. 
$C^1$ finite element spaces with a stable splitting are also required in Böhmer's proofs of the error bounds for his 
method of numerical solution of fully nonlinear elliptic equations \cite{Boehmer08}.
A stable splitting of the space $\mathbb{S}^{1,2}_{5,I}(\tri)$ will be obtained if a stable local basis for a stable complement of 
$\mathbb{S}^{1,2}_{5,0}(\tri)$ in $\mathbb{S}^{1,2}_{5,I}(\tri)$ is constructed, which we leave to a future work.

\section{Interpolation operator \label{operator}}\mprup{operator}

We denote by $\p^\alpha f$, $\alpha\in\mathbb{Z}^2_+$, the partial derivatives of $f$ and consider the
usual Sobolev spaces $H^m(\Omega)$ with the seminorm and norm defined by
$$
|f|_{H^m(\Omega)}^2=\sum_{|\alpha|=m}\|\p^\alpha f\|_{L^2(\Omega)}^2,\quad
\|f\|_{H^m(\Omega)}^2=\sum_{k=0}^m|f|_{H^k(\Omega)}^2\quad (H^0(\Omega)=L^2(\Omega)),$$
where $|\alpha|:=\alpha_1+\alpha_2$. We set
$H^1_0(\Omega)=\{f\in H^1(\Omega):f|_{\p\Omega}=0 \}$.

In this section we construct an interpolation operator $I_\tri:H^5(\Omega)\cap H^1_0(\Omega)\to \mathbb{S}^{1,2}_{5,0}(\tri)$
and estimate its error for the functions in
$H^m(\Omega)\cap H^1_0(\Omega)$, $m=5,6$, in the next section.

As in \cite{DKS}  we choose domains $\Omega_j\subset\Omega$, $j=1,\ldots,n$, with Lipschitz boundary such that
\begin{itemize}
\item[(a)]\quad $\p\Omega_j\cap\p\Omega = \Gamma_j$,
\item[(b)]\quad $\p\Omega_j\setminus\p\Omega$ is composed of a finite number of straight line segments,
\item[(c)]\quad $q_j(x)>0$ for all $x\in \overline{\Omega}_j\setminus \Gamma_j$, and
\item[(d)]\quad $\Omega_j\cap\Omega_k=\emptyset$ for all $j\ne k$.
\end{itemize}
In addition we assume that the triangulation $\tri$ is such that %
\begin{itemize}
\item[(e)]\quad $\overline{\Omega}_j$ contains every triangle $T\in\tri_P$ whose curved edge is part of
$\Gamma_j$,
\end{itemize}
and that $q_j$ satisfy (without loss of generality)
\begin{itemize}
\item[(f)]\quad $\displaystyle\max_{x\in\overline{\Omega}_j}\|\nabla q_j(x)\|_2\le 1\text{ and }
\|\nabla^2 q_j\|_2\le1,\;\text{for all }j =1,\ldots,n$,
\end{itemize}
where $\nabla^2 q_j$ denotes the (constant) Hessian matrix of $q_j$.

Note that (e) will hold with the same set $\{\Omega_j: j=1,\ldots,n\}$
for any triangulations obtained by subdividing the triangles
of $\tri$.

The following lemma can be shown following the lines of the proof of \cite[Theorem~6.1]{HRW01},
see also  \cite[Theorem 3.1]{DKS}.

\begin{lemma}%
\label{Hardyl}%
There is a constant $K$ depending only on $\Omega$, the choice of $\Omega_j$, $j=1,\ldots,n$, and $m\ge1$,
 such that for all $j$ and $u\in H^m(\Omega)\cap H^1_0(\Omega)$,
\begin{equation} \label{Hardye}
|u/q_j|_{H^{m-1}(\Omega_j)}\le K \|u\|_{H^{m}(\Omega_j)}.
\end{equation}
\end{lemma}

Given a a unit vector $\tau=(\tau_x,\tau_y)$ in the plane,
we denote by $D_\tau$ the directional derivative operator in the direction of
$\tau$ in the plane, so that
$$
D_\tau f:= \tau_xD_x f+\tau_y D_y f, \quad D_xf:=\p f/\p x,
 \quad D_yf:=\p {f}/\p{y}.
$$
Given $f\in C^{\alpha+\beta}(\tri),~\alpha,\beta\geq 0$, any number
$$
\eta f=D^\alpha_{\tau_1}D^\beta_{\tau_2}(f|_T)(z),
$$
where $T\in \tri,~z\in T$, and $\tau_1,\tau_2$ are some unit vectors in the plane,
is said to be a \emph{nodal value} of $f$, and the linear functional $\eta:C^{\alpha+\beta}(\tri)\rightarrow \mathbb{R}$
is a \emph{nodal functional}, with $d(\eta):=\alpha+\beta$ being the \emph{degree} of $\eta$.

For some special choices of $z,\tau_1,\tau_2$, we use the following notation:
\begin{itemize}
\item If $v$ is a vertex of $\tri$ and $e$ is an edge attached to $v$, we set
$$
D_e^\alpha f(v):= D_\tau^\alpha(f|_T)(v),\quad \alpha\geq 1,
$$
where $\tau$ is the unit vector in the direction of $e$ away from $v$, and $T\in \tri$ is one of the triangles with edge $e$.
\item If $v$ is a vertex of $\tri$ and $e_1,e_2$ are two consecutive edges attached to $v$, we set
$$
D^\alpha_{e_1}D^\beta_{e_2}f(v):=D^\alpha_{\tau_1}D^\beta_{\tau_2}(f|_T)(v),\quad \alpha,\beta\geq 1,
$$
where $T\in \tri$ is the triangle with vertex $v$ and edges $e_1,e_2$, and $\tau_i$ is the unit vector in the $e_i$ direction
away from $v$.
\item For every edge $e$ of the triangulation $\tri$ we choose a unit vector $\tau^\bot$ (one of two possible) orthogonal to
$e$ and set
$$
D^\alpha_{e^\bot}f(z):=D^\alpha_{\tau^\bot}f(z),\quad z\in e,\quad \alpha\geq 1,
$$
provided $f\in C^\alpha(z)$.
\end{itemize}

On every edge $e$ of $\tri$, with vertices $v'$ and $v''$, we define three points on $e$ by
$$
z_e^j:=v'+\frac{j}{4}(v''-v'),\quad j=1,2,3.
$$

\begin{figure}[htbp!]
\centering
\begin{tabular}{c}
\psfrag{v1}{$v_1$}
\psfrag{v2}{$v_2$}
\psfrag{v3}{$v_3$}
\includegraphics[height=0.3\textwidth]{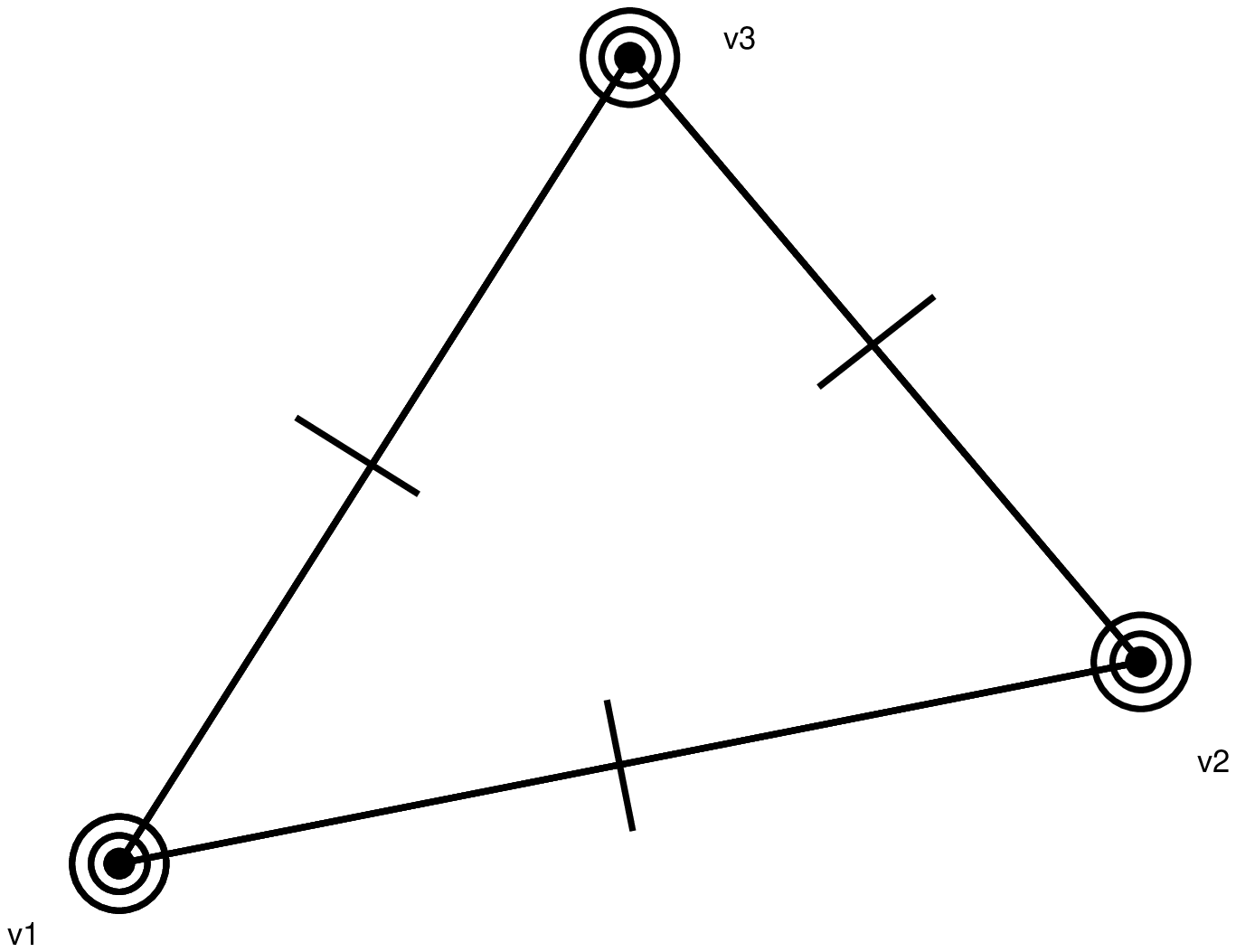}
\qquad
\psfrag{v1}{$v_1$}
\psfrag{v2}{$v_2$}
\psfrag{v3}{$v_3$}
\includegraphics[height=0.3\textwidth]{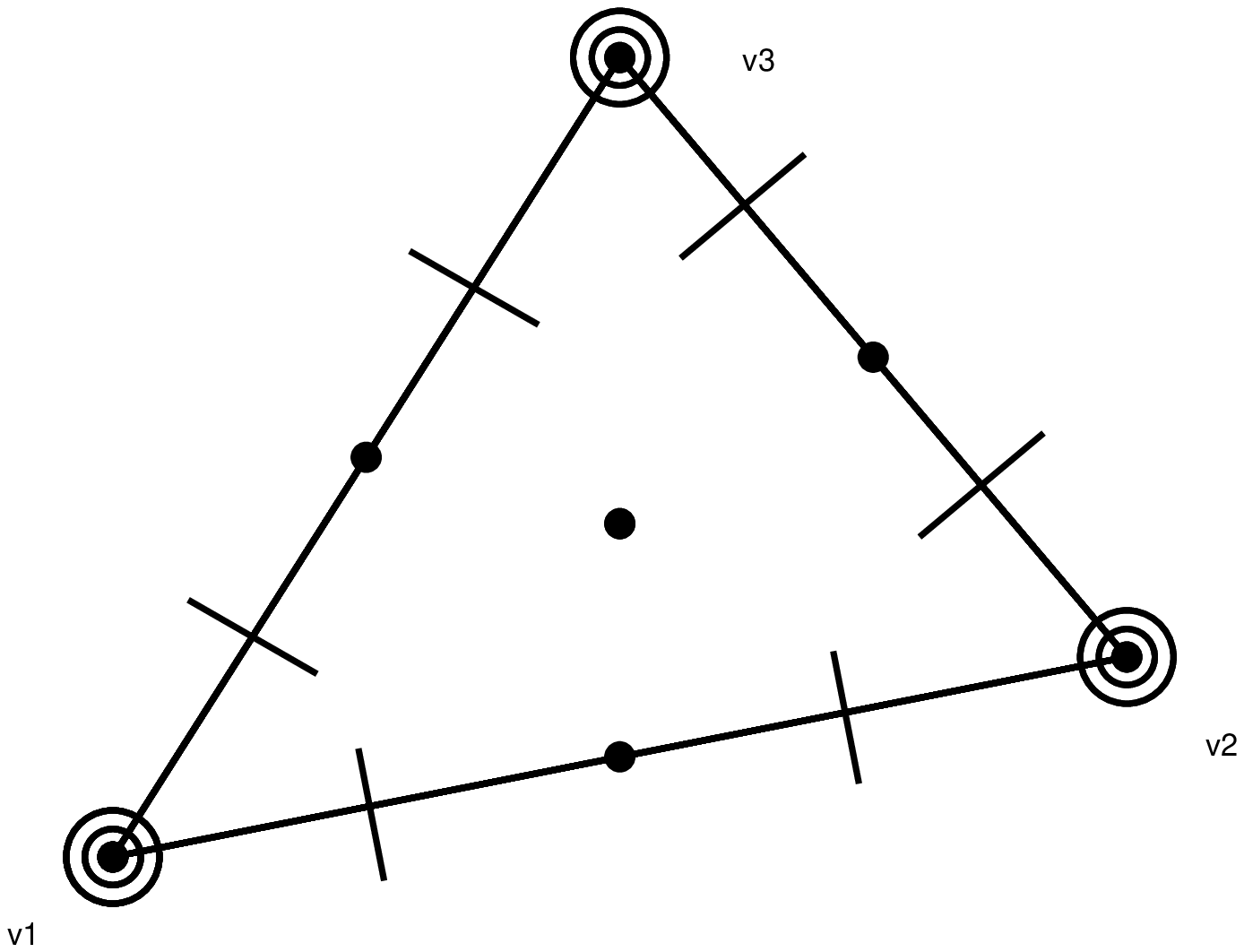}
\end{tabular}
\caption{Nodal functionals corresponding to $\CN_T^0$ (left) and $\CN_T^B$ (right).}
\label{ori}
\end{figure}

For every triangle $T\in \tri_0$ with vertices $v_1,v_2,v_3$ and edges $e_1,e_2,e_3$, we define $\CN_T^0$ to
be the set of nodal functionals corresponding to the nodal values
$$D_x^\alpha D_y^\beta f(v_i), \quad 0\leq \alpha+\beta \leq 2,\quad i=1,2,3,$$
$$
D_{e_i^\bot}f(z_{e_i}^2),\quad i=1,2,3,
$$
see Figure~\ref{ori} (left), where the nodal functionals are depicted in the usual way
by dots, segments and circles as for example in \cite{Ciarlet}.

\begin{figure}[htbp!]
\centering
\begin{tabular}{c}
\psfrag{v1}{$v_1$}
\psfrag{v2}{$v_2$}
\psfrag{v3}{$v_3$}
\includegraphics[height=0.35\textwidth]{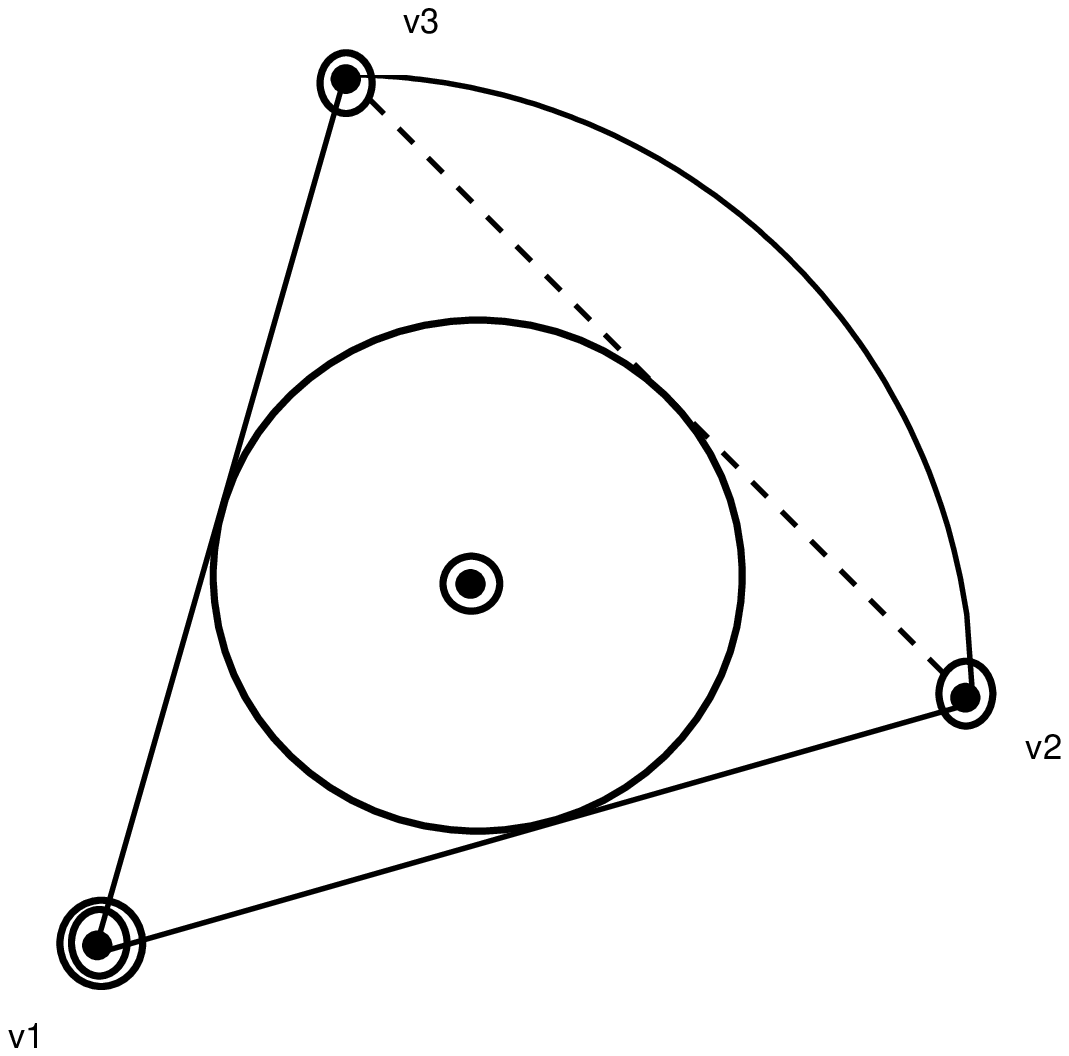}
\qquad
\psfrag{v1}{$v_1$}
\psfrag{v2}{$v_2$}
\psfrag{v3}{$v_3$}
\includegraphics[height=0.35\textwidth]{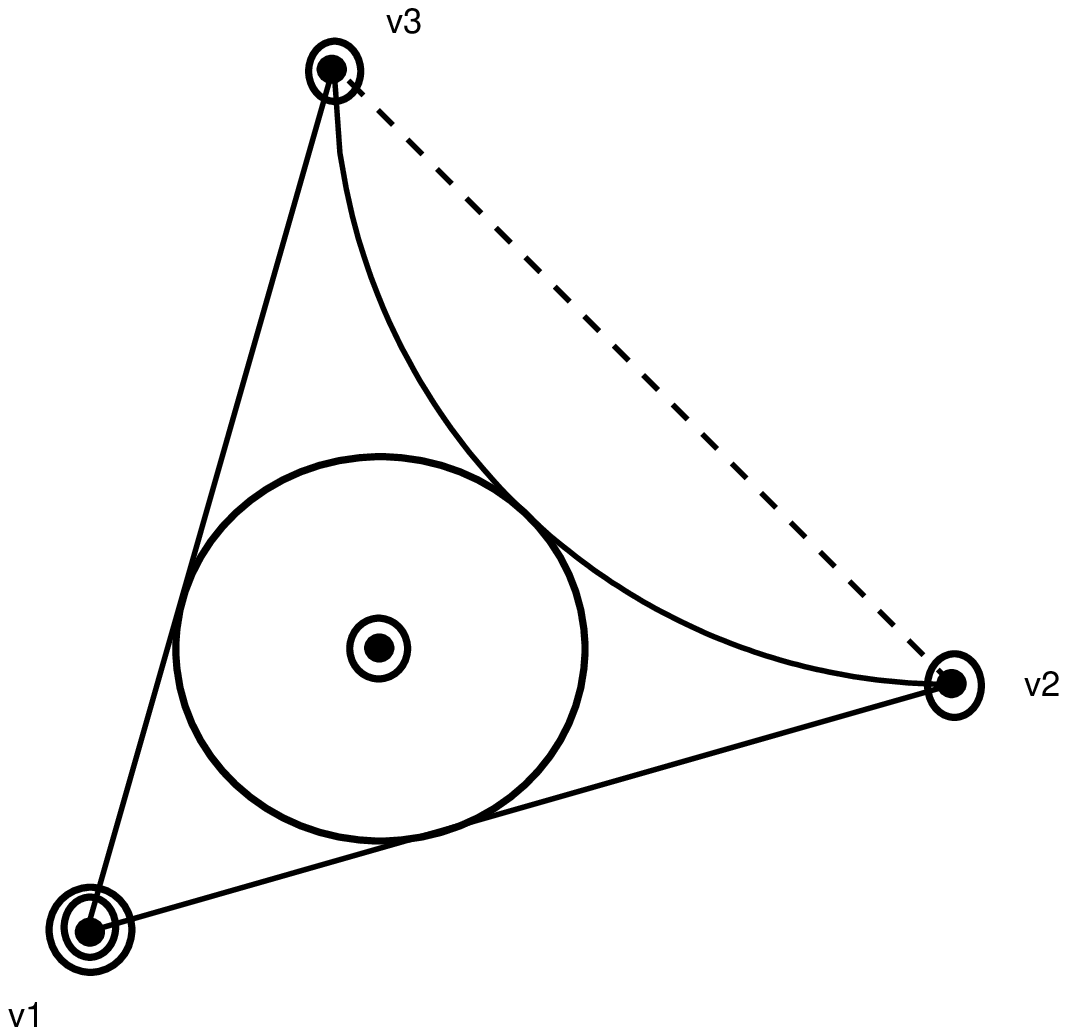}
\end{tabular}
\caption{Nodal functionals corresponding to $\CN_T^P$.}
\label{pief}
\end{figure}

Let $T\in \tri_P$.
We define $\CN_T^P$ to be the set of nodal functionals corresponding to the nodal values
$$D_x^\alpha D_y^\beta f(v_1), \quad 0\leq \alpha+\beta \leq 2,$$
$$
D^\alpha_{x}D^\beta_{y}f(v_i), \quad 0\leq \alpha+\beta \leq 1,\quad i=2,3,
$$
$$D_x^\alpha D_y^\beta f(c_T), \quad 0\leq \alpha+\beta \leq 1,$$
where $v_1$  the interior vertex of $T$,  $v_2,v_3$ are  boundary vertices, and
$c_T$ is the center of the disk $B_T$,
see Figure~\ref{pief}.

Let $T\in \tri_B$ with vertices $v_1,v_2,v_3$. %
We define $\CN_T^{B,1}$ to be the set of nodal functionals corresponding to the nodal
value
$$
f(c_T),\quad c_T:=(v_1+v_2+v_3)/3.$$
Also we define $\CN_T^{B,2}$ to be the set of nodal functionals corresponding to the nodal values
$$
f(z_{e_i}^2), \quad i=1,2,3,$$
$$
D_x^\alpha D_y^\beta f(v_i), \quad 0\leq \alpha+\beta \leq 2,\quad i=1,2,3,
$$
$$
D_{e_i^\bot}f(z_{e_i}^j), \quad j=1,3,\quad i=1,2,3,
$$
where $v_1$ is the boundary vertex and $v_2,v_3$ are the interior vertices of $T$.
We set
$$
\CN_T^B:=\CN_T^{B,1}\cup \CN_T^{B,2},$$
see Figure~\ref{ori} (right).

We define an operator %
$I_\tri: H^5(\Omega)\cap H^1_0(\Omega)\to \sps$ of interpolatory type.
Let $u\in H^5(\Omega)\cap H^1_0(\Omega)$. By Sobolev embedding we assume without loss of generality
that $u\in C^3(\overline{\Omega})$. For all $T\in\tri_0\cup\tri_P$ we set $I_\tri u|_T=I_T(u|_T)$, with the local operators
$I_T$ defined as follows.

If $T\in \tri_0$, then $p:=I_Tu$ is the polynomial of degree $5$ that satisfies the following
interpolation conditions:
$$
\eta p= \eta u, \quad \textrm{ for all }\eta\in \CN_T^0.
$$
This is a well-known Argyris interpolation scheme, see e.g.~\cite[Section 6.1]{LSbook}, which
guarantees the existence and uniqueness of the polynomial $p$.

Let $T\in\tri_P$ with the curved edge on $\Gamma_j$.
Then $I_Tu:=pq_j$, where $p\in\mathbb{P}_{4}$ satisfies the following interpolation condition:
\begin{equation}\label{intpie}
\eta p =\eta (u/q_j), \quad \textrm{ for all }\eta\in \CN_T^P.
\end{equation}
The nodal functionals in $\CN_T^P$ are well defined for $u/q_j$
even though the vertices $v_2,v_3$ of $T$
lie on the boundary $\Gamma_j$ because $u/q_j\in H^4(\Omega_j)$ by Lemma~\ref{Hardyl} and hence
$u/q_j$ may be identified with a function $\tilde u\in C^2(\overline{\Omega}_j)$ by Sobolev embedding.
The interpolation scheme \eqref{intpie} defines a unique polynomial $p\in \mathbb{P}_{4}$, which
will be justified  in the proof of Lemma~\ref{pieb0}.
In addition, we will need the following statement.
\begin{lemma} \label{piel} The polynomial $p$ defined by \eqref{intpie} satisfies
$$ %
D^\alpha_xD^\beta_y(pq_j)(v)=D^\alpha_xD^\beta_yu(v),\quad 0\le \alpha+\beta\le2,
$$ %
where $v$ is any vertex of the pie-shaped triangle $T$.
\end{lemma}
\pf
By \eqref{intpie}, $p(v)q_j(v)=\tilde u(v)q_j(v)=u(v)$, where $\tilde u\in C^2(\overline{\Omega}_j)$
is the above function satisfying $u=\tilde uq_j$. Moreover,
\begin{align*}
\nabla (pq_j)(v)&=\nabla p(v)q_j(v)+p(v)\nabla q_j(v)\\
&=\nabla \tilde u(v)q_j(v)+\tilde u(v)\nabla q_j(v)\\
&=\nabla (\tilde uq_j)(v) = \nabla u(v).
\end{align*}
Similarly, if $v$ is the interior vertex of $T$, then
\begin{align*}
\nabla^2 (pq_j)(v)&=\nabla^2 p(v)q_j(v)+\nabla p(v)(\nabla q_j(v))^T+\nabla q_j(v)(\nabla p(v))^T+p(v)\nabla^2 q_j(v)\\
&=\nabla^2 \tilde u(v)q_j(v)+\nabla \tilde u(v)(\nabla q_j(v))^T+\nabla q_j(v)(\nabla \tilde u(v))^T+\tilde u(v)\nabla^2 q_j(v)\\
&=\nabla^2 u(v).
\end{align*}
If $v$ is one of the boundary vertices, then $q_j(v)=0$, and hence
\begin{align*}
\nabla^2 (pq_j)(v)&=\nabla p(v)(\nabla q_j(v))^T+\nabla q_j(v)(\nabla p(v))^T+p(v)\nabla^2 q_j(v)\\
&=\nabla \tilde u(v)(\nabla q_j(v))^T+\nabla q_j(v)(\nabla \tilde u(v))^T+\tilde u(v)\nabla^2 q_j(v)\\
&=\nabla^2 u(v).\meop
\end{align*}
It is easy to deduce from Lemma~\ref{piel} that the interpolation conditions
for $p$ at the boundary vertices $v_2,v_3$ of $T$ can be equivalently formulated as follows: For $i=2,3$,
\begin{align}\label{intp}
\begin{split}
p(v_i)&=\frac{\p u}{\p n_i}(v_i)\Big/\frac{\p q_j}{\p n_i}(v_i),\\
\frac{\p p}{\p n_i} (v_i)=\frac{1}{2}\frac{\p^2 u}{\p n_i^2}(v_i)\Big/\frac{\p q_j}{\p n_i}&(v_i),\qquad
\frac{\p p}{\p \tau_i} (v_i)=\frac{\p^2 u}{\p n_i\p \tau}(v_i)\Big/\frac{\p q_j}{\p n_i}(v_i),
\end{split}
\end{align}
where $n_i$ and $\tau_i$ are the normal and the tangent unit vectors to the curve $q_j(x)=0$ at $v_i$.

Finally, assume that $T\in\tri_B$ with vertices $v_1,v_2,v_3$ where $v_1$ is a boundary vertex.
Then $I_Tu=p\in\mathbb{P}_{6}$ satisfies the following interpolation conditions:
$$
\eta p =\eta u, \quad \textrm{ for all }\eta \in \CN_T^{B,1},
$$
and
$$
\eta p =\eta I_{T_i}u, \quad \textrm{ for all }\eta \in \CN_i \subset \CN_T^{B,2}, \quad i=1,2,3,
$$
where $T_1$ is a triangle in $\tri_0$ sharing an edge $e_1=\langle v_2, v_3 \rangle$ with $T$ and $\CN_1$ corresponds to the nodal values
$$
f(z_{e_1}^2),\quad D_{e_1^\bot}f(z_{e_1}^i), \quad i=1,3,
$$
$$D_x^\alpha D_y^\beta f(v_i), \quad 0\leq \alpha+\beta \leq 2,\quad i=2,3;$$
$T_2$ is a triangle in $\tri_P$ sharing an edge $e_2=\langle v_1, v_2 \rangle$ with $T$ and
$\CN_2$ corresponds to the nodal values
$$
f(z_{e_2}^2),\quad D_{e_2^\bot}f(z_{e_2}^i), \quad i=1,3,
$$
$$
D^\alpha_{x}D^\beta_{y}f(v_1), \quad 0\leq \alpha+\beta \leq 2;
$$
and $T_3$ is a triangle in $\tri_P$ sharing an edge $e_3=\langle v_1, v_3 \rangle$ with $T$ and
$\CN_3$ corresponds to the nodal values
$$
f(z_{e_3}^2),\quad D_{e_3^\bot}f(z_{e_3}^i), \quad i=1,3.
$$
Since $\CN_T^{B,2}=\CN_1\cup \CN_2\cup \CN_3$ and $\CN_T^{B}=\CN_T^{B,1}\cup\CN_T^{B,2}$ is a well posed interpolation scheme \cite{Sch89}
for polynomials of degree 6, it follows that $p$ is uniquely defined by the above conditions.

\begin{theorem} \label{C1}
Let $u\in H^5(\Omega)\cap H^1_0(\Omega)$. Then  $I_\tri u\in\sps$.
\end{theorem}

\pf
By the above construction $I_\tri u$ is a piecewise polynomial of degree 5 on all triangles in $\tri_0$
and degree 6 on the triangles in $\tri_P\cup\tri_B$. Moreover, $I_\tri u$ vanishes on the boundary of $\Omega$.

To see that $I_\tri u\in\sps$ we thus need to show the $C^1$ continuity of $I_\tri u$ across all interior edges of $\tri$. If $e$ is a
common edge of two triangles $T',T^{''}\in \tri_0$, then the $C^1$ continuity follows from the standard argument for $C^1$ Argyris
finite element, see \cite[Chapter 3]{BrennerScott} and \cite[Section 6.1]{LSbook}.

Next we will show the $C^1$ continuity of $I_\tri u$ across edges shared by buffer triangles with either ordinary or pie-shaped
triangles. Let $T\in \tri_B$ and
$T'\in \tri_0\cup \tri_P$ with common edge $e'=\langle v',v''\rangle$,
and let $p=I_{T}u$ and $s=I_{T'}u$. Consider the
univariate polynomials $p|_{e'}$ and $s|_{e'}$ and let $q=p|_{e'}-s|_{e'}$. Assuming that the edge $e'$ is
parameterized by $t\in[0,1]$,
Then $q$ is a univariate polynomial of degree 6 with respect to the parameterization $v'+t(v''-v')$,
$t\in[0,1]$. Similarly, we consider the
orthogonal/normal derivatives $D_{e'^\bot}p|_{e'}$ and $D_{e'^\bot}s|_{e'}$ and let
$r=D_{e'^\bot}p|_{e'}-D_{e'^\bot}s|_{e'}$, then $r$ is a univariate polynomial of degree 5 with respect to the
same parameter $t$.
The $C^1$ continuity will follow if we show that both $q$ and $r$ are zero functions.

If $T'=T_1\in \tri_0$, then
using the interpolation conditions corresponding to $\CN_1\subset\CN_T^{B,2}$, we have
\begin{align*}
&q(0)=q'(0)=q''(0)=q(1/2)=q(1)=q'(1)=q''(1)=0,\\
&r(0)=r'(0)=r(1/4)=r(3/4)=r(1)=r'(1)=0,
\end{align*}
which implies $q\equiv 0$ and $r\equiv 0$.

If $T'=T_2\in \tri_P$, then the interpolation conditions corresponding to $\CN_2\subset\CN_T^{B,2}$ imply
\begin{align*}
&q(0)=q'(0)=q''(0)=q(1/2)=0,\\
&r(0)=r'(0)=r(1/4)=r(3/4)=0,
\end{align*}
In view of Lemma~\ref{piel}, we have
$$
D^\alpha_xD^\beta_ys(v_2)=D^\alpha_xD^\beta_yu(v_2)
=D^\alpha_xD^\beta_yp(v_2),\quad 0\le \alpha+\beta\le2,$$
which implies
$$
q(1)=q'(1)=q''(1)=0,\quad r(1)=r'(1)=0,$$
and hence $q\equiv 0$ and $r\equiv 0$.

If $T'=T_3\in \tri_P$, then the interpolation conditions corresponding to $\CN_3\subset\CN_T^{B,2}$
imply
$$
q(1/2)=0,\quad r(1/4)=r(3/4)=0,$$
whereas Lemma~\ref{piel} gives
\begin{align*}
&q(0)=q'(0)=q''(0)=0,\quad r(0)=r'(0)=0,\\
&q(1)=q'(1)=q''(1)=0,\quad r(1)=r'(1)=0,
\end{align*}
which completes the proof.\eop

In follows from Lemma~\ref{piel} that %
$I_\tri u$ is twice differentiable at the boundary vertices, and thus
$$
I_\tri u\in\{s\in  \mathbb{S}^1_5(\tri)\sp s \text{ is twice differentiable at any vertex and } s|_\Gamma = 0\}.$$
Moreover, $I_\tri u$  satisfies the following interpolation conditions:
 $$
D_x^\alpha D_y^\beta I_\tri u(v)=D_x^\alpha D_y^\beta u(v), \quad 0\leq \alpha+\beta \leq 2,\quad
\text{for all $v\in V$},$$
$$
D_{e^\bot}I_\tri u(z_{e}^2)=D_{e^\bot}u(z_{e}^2), \quad \text{for all edges $e$ of $\tri_0$},
$$
$$
D_x^\alpha D_y^\beta I_\tri u(c_T)=D_x^\alpha D_y^\beta u(c_T), \quad 0\leq \alpha+\beta \leq 1,\quad
\text{for all $T\in \tri_P$},$$
$$
I_\tri u(c_T)=u(c_T), \quad 
\text{for all $T\in \tri_B$},$$
where $c_T$ denotes the center of the disk $B_T$ inscribed into $T^*$ if $T$ is a pie-shaped triangle, and the
barycenter of $T$ if $T$ is a buffer triangle. In view of \eqref{intp}, $I_\tri u\in\sps$ is uniquely defined by
these conditions for any $u\in C^2(\overline{\Omega})$.

\section{Error bounds \label{bounds}}\mprup{bounds}

In this section we estimate the error $\|u-I_\tri u\|_{H^k(\Omega)}$
for functions $u\in H^m(\Omega)\cap H^1_0(\Omega)$, $m=5,6$. Similar to \cite[Section 3]{DKS}, 
we follow the standard finite element techniques
involving the Bramble-Hilbert Lemma  (see \cite[Chapter 4]{BrennerScott}) on the ordinary triangles,
and make use of  the estimate \eqref{Hardye} on the pie-shaped triangles. 
Since the spline $I_\tri u$ on the buffer triangles is constructed in part by interpolation and in part by
the smoothness conditions, the estimate 
of the error on such triangles relies in particular
on the estimates of the interpolation error  on the 
neighboring ordinary and buffer triangles.

\begin{lemma} \label{pieb0} If $p\in\mathbb{P}_4$ and $T \in \tri_P$, then
\begin{equation}\label{intb0}
\|p|_{T^{*}}\|_{L^\infty(T^{*})}\le \max_{\eta\in \CN_T^P}h_{T^{*}}^{d(\eta)}|\eta p|,
\end{equation}
where $T^*$ is the triangle obtained by replacing the curved edge of $T$ by the straight line segment,
 and $h_{T^{*}}$ is the diameter of $T^{*}$.
Similarly, if $p\in\mathbb{P}_6$ and $T \in \tri_B$, then
\begin{equation}\label{intb1}
\|p|_T\|_{L^\infty(T)}\le \max_{\eta\in \CN_T^B}h_{T}^{d(\eta)}|\eta p|,
\end{equation}
where $h_{T}$ is the diameter of $T$.
\end{lemma}
\pf %
To show the estimate \eqref{intb0} for $T^{*}$, we follow the  proof of \cite[Lemma
3.9]{DNZ}. We note that we only need to show that the interpolation scheme for pie-shaped triangles is a valid
scheme, that is, we need to show that $\CN_T^P$ is $\mathbb{P}_{4}$-{unisolvent}, and the rest of the proof
can be done similar to that of \cite[Lemma 3.9]{DNZ}. 
Recall that a set of functionals $\CN$ is said to be $\mathbb{P}_{d}$-\emph{unisolvent} 
if the only polynomial $p\in \mathbb{P}_{d}$ satisfying $\eta p=0$ for 
$\eta \in \CN$ is the zero
function.

Let $T^{*}=\langle v_1,v_2,v_3\rangle$, where $v_1$ is the interior vertex. Set $e_1:=\langle v_1,v_2
\rangle$,  $e_2:=\langle v_2,v_3\rangle$, $e_3:=\langle v_3,v_1\rangle$, see Figure~\ref{pief}. The
interpolation conditions along $e_1, e_3$ imply that $s$ vanishes on these edges. After splitting out the
linear polynomials factors which vanish along the edges $e_1, e_3$ we obtain a valid interpolation scheme for
quadratic polynomials with function values at the three vertices, and function and gradient values at the
the barycenter $c$ of $B_T\subset T^*$. The validity of this scheme can be seen by looking at a straight line $\ell$
through $c$ and any one of the vertices of $T^{*}$. Along the line $\ell$, a function value
is given at the vertex and a function value together with the first derivative are given at the point $c$, and this
set of data is {unisolvent} for the univariate quadratic polynomials, which means $s$ must vanish along
$\ell$. After factoring out the  respective linear polynomial, we are left with function values at three
non-collinear points, which defines a valid interpolation scheme for the remaining linear polynomial factor
of $s$.

To show the estimate \eqref{intb1} for $T \in \tri_B$, the proof is similar. We need
to show the set of functionals $\CN_T^B$ is $\mathbb{P}_{6}$-{unisolvent} but this follows from the standard scheme of
\cite{Sch89} for polynomials of degree six.

We note that the argument of the proof of \cite[Lemma 3.9]{DNZ} applies to  affine invariant 
interpolation schemes, that is the schemes that use
the edge derivatives. As our scheme relies on the standard derivatives 
in the direction of the $x,y$ axes, we need to express the 
edge derivatives as linear combinations of the $x,y$ derivatives as follows.
Assume that $e_1,e_2$ are two edges that emanate from a vertex $v$. Let $\tau_i=(\tau_{i1},\tau_{i2})$
be the unit vector in the direction of $e_i$ away from $v$, $i=1,2$.
Then we can easily obtain the following identities
$$D_{e_i}f(v)=\tau_{i1}D_xf(v)+\tau_{i2}D_yf(v),$$
$$
D^2_{e_i}f(v)=\tau_{i1}^2 D^2_xf(v)+2\tau_{i1}\tau_{i2}D_xD_yf(v)+\tau_{i2}^2D^2_yf(v),
$$
$$
D_{e_1}D_{e_2}f(v)=\tau_{11}\tau_{21}D_x^2f(v)+(\tau_{11}\tau_{22}+\tau_{12}\tau_{21})D_xD_yf(v)+\tau_{12}\tau_{22}D_y^2f(v).
\meop$$

\begin{lemma} \label{pie}
Let $T\in\tri_P$ and its curved edge $e\subset\Gamma_j$. Then
\begin{equation}\label{pieb1}
\|I_Tu\|_{L^\infty(T)}\le C_1\max_{0\le \ell\le 2}h_{T}^{\ell+1}|u/q_j|_{W^{\ell}_\infty(T)}
\quad\text{if }\;u\in H^5(\Omega)\cap H^1_0(\Omega),
\end{equation}
where $C_1$ depends only on $h_T/\rho_T$.
Moreover, if $5\le m\le 6$, %
then for any $u\in H^m(\Omega)\cap H^1_0(\Omega)$,
\begin{align}\label{pieb2}
\|u-I_{T}u\|_{H^k(T)}&\le C_2h_T^{m-k}|u/q_j|_{H^{m-1}(T)},\quad k=0,\ldots,m-1,\\
\label{pieb3}
|u-I_{T}u|_{W^k_\infty(T)}&\le C_3h_T^{m-k-1}|u/q_j|_{H^{m-1}(T)},\quad k=0,\ldots,m-2,
\end{align}
where $C_2,C_3$ depend only on $h_T/\rho_T$.
\end{lemma}

\pf We will denote by $\tilde C$ constants which may depend only on $h_T/\rho_T$ and on
$\Omega$.
Assume that $u\in H^5(\Omega)\cap H^1_0(\Omega)$ and recall that by definition $I_Tu=pq_j$,
where $p\in\mathbb{P}_{4}$ satisfies the interpolation conditions \eqref{intpie}. %
Since $u\in H^5(\Omega_j)\cap H^1_0(\Omega_j)$, it follows that $u/q_j\in H^4(\Omega_j)$ by
Lemma~\ref{Hardyl}, and hence
$u/q_j \in C^2(\overline{\Omega}_j)$ by Sobolev embedding.
From Lemma \ref{pieb0} we have
\begin{equation}\label{intb}
\|p\|_{L^\infty(T^{*})}\le \max_{\eta\in \CN_T^P}h_{T^{*}}^{d(\eta)}|\eta p|,
\end{equation}
and hence
\begin{align*}
\|p\|_{L^\infty(T^{*})}\le \max_{\eta\in \CN_T^P}h_{T^{*}}^{d(\eta)}|\eta (u/q_j)|
\le \tilde C \max_{0\le \ell\le 2}h_{T}^{\ell}|u/q_j|_{W^{\ell}_\infty(T)}.
\end{align*}
As in the proof of \cite[Theorem 3.2]{DKS}, we can show that for any polynomial of degree at most 6,
\begin{equation}\label{normTT}
\|s\|_{L^\infty(T)}\le\tilde C\|s\|_{L^\infty(T^{*})}\quad\text{and}\quad\|s\|_{L^\infty(T^{*})}\le\tilde C\|s\|_{L^\infty(T)}.
\end{equation}
By using (f) it is easy to show that $\|q_j\|_{L^\infty(T)}\le h_T$, and hence
$$
\|I_Tu\|_{L^\infty(T)}=\|pq_j\|_{L^\infty(T)}\le h_T\|p\|_{L^\infty(T)},$$
which completes the proof of \eqref{pieb1}.

Moreover, since the area of $T$ is less than or equal $\frac{\pi}{4}h^2_T$ and $\p^\alpha (I_Tu)\in\mathbb{P}_{6-k}$
if $|\alpha|=k$, it follows that
$$
\|\p^\alpha (I_Tu)\|_{L^2(T)}\le \frac{\sqrt{\pi}}{2}h_T\|\p^\alpha (I_Tu)\|_{L^\infty(T)}
\le \tilde Ch_T\|\p^\alpha (I_Tu)\|_{L^\infty(T^{*})}.$$
By Markov inequality (see e.g.~\cite[Theorem 1.2]{LSbook}) we get furthermore
$$
\|\p^\alpha (I_Tu)\|_{L^\infty(T^{*})}\le
\tilde C \rho_T^{-k}\|I_Tu\|_{L^\infty(T^{*})},$$
and hence in view of \eqref{normTT}
$$ %
|I_Tu|_{H^k(T)}\le \tilde C h^{1-k}_T\|I_Tu\|_{L^\infty(T)}.
$$ %
In view of  \eqref{pieb1} we  arrive at
\begin{equation}\label{pieb4}
|I_Tu|_{H^k(T)}\le \tilde C\max_{0\le \ell\le 2}h_{T}^{\ell+2-k}|u/q_j|_{W^{\ell}_\infty(T)},\quad \text{if }\;u\in
H^5(\Omega)\cap H^1_0(\Omega).
\end{equation}

Let $m\in\{5,6\}$, and let $u\in H^m(\Omega)\cap H^1_0(\Omega)$. It follows from
Lemma~\ref{Hardyl} that $u/q_j\in H^{m-1}(T)$.
By the results in \cite[Chapter 4]{BrennerScott} there exists a polynomial $\tilde p\in\mathbb{P}_{m-2}$ such that
\begin{align}
\begin{split}\label{BH}
\|u/q_j-\tilde p\|_{H^k(T)}&\le \tilde Ch_T^{m-k-1}|u/q_j|_{H^{m-1}(T)},\quad k=0,\ldots,m-1,\\
|u/q_j-\tilde p|_{W^k_\infty(T)}&\le \tilde Ch_T^{m-k-2}|u/q_j|_{H^{m-1}(T)},\quad k=0,\ldots,m-2.
\end{split}
\end{align}
Indeed, a suitable $\tilde p$ is given by the
\emph{averaged Taylor polynomial} \cite[Definition 4.1.3]{BrennerScott} with respect to the disk $B_T$,
and the inequalities in \eqref{BH} follow from \cite[Lemma 4.3.8]{BrennerScott} (Bramble-Hilbert Lemma) and
an obvious generalization of
\cite[Proposition 4.3.2]{BrennerScott}, respectively. It is easy to check by inspecting the proofs in
\cite{BrennerScott} that the quotient $h_T/\rho_T$ can be used in the estimates instead of the
chunkiness parameter used there.

Since
$$
u-I_Tu=(u/q_j-\tilde p)q_j-I_T(u - \tilde pq_j), $$
we have for any norm $\|\cdot\|$,
$$
\|u-I_Tu\|\le\|(u/q_j-\tilde p)q_j\|+\|I_T(u - \tilde pq_j)\|. $$
In view of (f) and \eqref{BH}, for any $k=0,\ldots,m-2$,
\begin{align*}
|(u/q_j-\tilde p)q_j|_{W^k_\infty(T)}&\le h_T |u/q_j-\tilde p|_{W^k_\infty(T)}+\|u/q_j-\tilde p\|_{W^{k-1}_\infty(T)}\\
&\le \tilde Ch_T^{m-k-1}|u/q_j|_{H^{m-1}(T)},
\end{align*}
and for any $k=0,\ldots,m-1$,
\begin{align*}
\|(u/q_j-\tilde p)q_j\|_{H^k(T)}
&\le \tilde Ch_T\|u/q_j-\tilde p\|_{H^k(T)}+\tilde C\|u/q_j-\tilde p\|_{H^{k-1}(T)}\\
&\le \tilde Ch_T^{m-k}|u/q_j|_{H^{m-1}(T)}.
\end{align*}

Furthermore, by the Markov inequality, \eqref{pieb1}, \eqref{pieb4} and \eqref{BH},
$$
|I_T(u - \tilde pq_j)|_{W^k_\infty(T)}\le \tilde C \max_{0\le \ell\le 2}h_{T}^{\ell+1-k}|u/q_j - \tilde p|_{W^{\ell}_\infty(T)}
\le \tilde Ch_T^{m-k-1}|u/q_j|_{H^{m-1}(T)},$$
$$
\|I_T(u - \tilde pq_j)\|_{H^k(T)}\le \tilde C\max_{0\le \ell\le 2}h_{T}^{\ell+2-k}|u/q_j - \tilde p|_{W^{\ell}_\infty(T)}
\le \tilde Ch^{m-k}_T|u/q_j|_{H^{m-1}(T)}.$$
By combining the inequalities in the five last displays we deduce \eqref{pieb2} and \eqref{pieb3}.
\eop

We are ready to formulate and prove our main result.

\begin{theorem} \label{approx}\mpr{aprox}
Let $5\le m\le 6$. For any $u\in H^m(\Omega)\cap H^1_0(\Omega)$,
\begin{equation}\label{appr1}
\Big(\sum_{T\in\tri}\|u-I_\tri u\|_{H^k(T)}^2\Big)^{1/2}\le
   Ch^{m-k}\|u\|_{H^m(\Omega)},\quad k=0,\ldots,m-1,
\end{equation}
where $h$ is the maximum diameter of the triangles in $\tri$,
and $C$ is a constant depending only on $\Omega$, the choice of $\Omega_j$, and the shape regularity constant
$R$ of $\tri$.
\end{theorem}

\pf
We estimate
the norms of $u-I_Tu$ on all triangles $T\in\tri$. The letter $C$ stands below for various constants depending only on the parameters mentioned
in the formulation of the theorem.

If $T\in\tri_0$, then $s|_T$ is a macro element as defined in \cite[Chapter 6]{LSbook}. Furthermore, by \cite[Theorem 6.3]{LSbook} the
set of linear functionals $\CN_T^0$ give rise to a stable local nodal basis, which is in particular uniformly
bounded. Hence by \cite[Theorem 2]{DY14} we obtain a Jackson estimate in the form
\begin{equation}\label{bound0}
\|u-I_{T}u\|_{H^k(T)}\le Ch_T^{m-k}|u|_{H^m(T)},\quad k=0,\ldots,m,
\end{equation}
where $C$ depends only on $h_T/\rho_T$.
If $T\in\tri_P$, with the curved edge $e\subset\Gamma_j$, %
then the Jackson estimate \eqref{pieb2}  holds by Lemma~\ref{pie}.

Let $T\in\tri_B$, $p:=I_\tri u|_T$ and let $\tilde p\in\mathbb{P}_6$ be the interpolation polynomial that satisfies
$\eta\tilde p=\eta u$ for all $\eta\in \CN_T^B $. Then
$$
\eta (\tilde p- p)=
\begin{cases} 0 & \text{if } \eta\in \CN_T^{B,1},\\
\eta (u-I_{T'}u) & \text{if } \eta\in \CN_T^{B,2},
\end{cases}$$
where  $T'=T'_\eta\in \tri_0\cup\tri_P$. Hence, by Markov inequality and \eqref{intb1} of Lemma \ref{pieb0},
we conclude that for $k=0,\ldots,m$,
$$
\|\tilde p-p\|_{H^k(T)}\le Ch_T^{1-k}\|\tilde p-p\|_{L^\infty(T)},$$
with
$$
\|\tilde p-p\|_{L^\infty(T)} \le
C\max\{h_T^{\ell}|u-I_{T'}u|_{W^\ell_\infty(T')}:0\le\ell\le2,\;T'\in \tri_0\cup\tri_P,\; T'\cap T\ne\emptyset\},
$$
whereas by the same arguments leading to \eqref{bound0} we have
$$
\|u-\tilde p\|_{H^k(T)}\le Ch_T^{m-k}|u|_{H^m(T)},$$
with the constants depending only on $h_T/\rho_T$.
If $T'\in\tri_0\cup\tri_P$, then by \eqref{pieb3} and the analogous estimate for $T'\in\tri_0$, compare
\cite[Corollary 4.4.7]{BrennerScott}, we have  for $\ell=0,1,2$,
$$
|u-I_{T'}u|_{W^\ell_\infty(T')}\le Ch_{T'}^{m-\ell-1}
\begin{cases}
|u|_{H^{m}(T')}&\text{if }T'\in\tri_0,\\
|u/q_j|_{H^{m-1}(T')}&\text{if }T'\in\tri_P,
\end{cases}$$
where $C$ depends only on $h_{T'}/\rho_{T'}$.
By combining these inequalities we obtain
an estimate of $\|u- I_Tu\|_{H^k(T)}$ by $C\tilde h^{m-k}$ times the maximum of $|u|_{H^{m}(T)}$,
$|u|_{H^{m}(T')}$ for $T'\in\tri_0$ sharing edges with $T$, and $|u/q_j|_{H^{m-1}(T')}$
for $T'\in\tri_P$ sharing edges with $T$. Here $C$ depends only on the maximum of
$h_{T}/\rho_{T}$ and $h_{T'}/\rho_{T'}$, and $\tilde h$ is the maximum of $h_T$ and all $h_{T'}$ for
$T'\in\tri_0\cup\tri_P$  sharing edges with $T$.

By using \eqref{bound0} on $T\in\tri_0$, \eqref{pieb2} on $T\in\tri_P$ and the estimate of the last
paragraph on $T\in\tri_B$, we get
\begin{align*}
\sum_{T\in\tri}\|u-I_\tri u\|_{H^k(T)}^2&\le Ch^{2(m-k)}\Big(\sum_{T\in \tri_0\cup \tri_B}|u|_{H^{m}(T)}^2
+\sum_{T\in\tri_P}|u/q_{j(T)}|_{H^{m-1}(T)}^2\Big),
\end{align*}
where $j(T)$ is the index of $\Gamma_j$ containing the curved edge
of $T\in\tri_P$. Clearly,
$$
\sum_{T\in\tri_0\cup \tri_B}|u|_{H^{m}(T)}^2\le|u|_{H^{m}(\Omega)}^2\le\|u\|_{H^{m}(\Omega)}^2,
$$
whereas by Lemma~\ref{Hardyl},
$$
\sum_{T\in\tri_P}|u/q_{j(T)}|_{H^{m-1}(T)}^2\le\sum_{j=1}^n|u/q_{j}|_{H^{m-1}(\Omega_j)}^2
\le K \|u\|_{H^m(\Omega)}^2,$$
where $K$ is the constant of \eqref{Hardye} depending only on $\Omega$ and the choice of $\Omega_j$.
\eop

\subsubsection*{Acknowledgements}
This research has been supported in part by the grant UBD/PNC2/2/RG/1(301) from Universiti Brunei
Darussalam.

\def\JAT{J.~Approx.\ Theory}
\def\JCAM{J.~Comput.\ Appl.\ Math.}
\def\AiCM{Advances in Comp.\ Math.}
\def\CA{Constr.\ Approx.}
\def\SJNA{SIAM J. Numer.\ Anal.}
\def\CAGD{CAGD}

\end{document}